\documentclass[11pt,a4paper]{article}
\usepackage{amsmath}
\usepackage{amsfonts}
 \usepackage{amsthm}
 \usepackage{hyperref}
 \usepackage[usenames]{color} 
 \usepackage{mathrsfs}
\theoremstyle{plain}
\newtheorem{thm}{Theorem}[section]
\newtheorem{lem}[thm]{Lemma}
\newtheorem{prop}[thm]{Proposition}

\theoremstyle{definition}
\newtheorem{rem}[thm]{Remark}

\renewcommand{\proof}{\vskip 0pt\noindent\textbf{\textit{Proof. }}}
\newcommand{\proofof}[1]{\vskip 0pt\noindent\textbf{\textit{Proof #1. }}}

\newcommand{\R}{\mathbb R}
\newcommand{\N}{\mathbb N}
\newcommand{\Z}{\mathbb Z}

\newcommand{\x}{{\bf x}}
\newcommand{\kk}{{\bf k}}

\renewcommand{\>}{\right>}
\newcommand{\abs}[1]{\left\vert#1\right\vert} 

\newcommand{\ind}{1\mkern-7mu1}

\numberwithin{equation}{section}
\newcommand{\La}{{\Lambda}}

\newcommand{\iy}{{\infty}}

\renewcommand{\P}{\mathbb{P}}
\newcommand{\pic}{{\pi^0_c}}
\newcommand{\piEc}{{\pi^{1,E}_{1-c}}}

\newcommand{\hpic}{{\hat \pi^0_c}}
\newcommand{\hpiEc}{{\hat \pi^{1,E}_{1-c}}}

 \title{\normalsize{\textbf{Diffusion limits at small times for coalescents with a  Kingman component}}}
 
  \author{
  \small{VLADA LIMIC}\thanks{Research supported in part by the ANR MANEGE grant.}\\
 \normalsize{CNRS UMR 8628, Laboratoire}\\
\normalsize{ de Math\'ematiques, Universit\'e }\\
 \normalsize{Paris-Sud, B\^atiment 425}\\
   \normalsize{91405 Orsay, France}\\
   \\
    \normalsize{e-mail: vlada.limic@math.u-psud.fr}
   \and
 \small{ANNA TALARCZYK}\thanks{Research supported in part by NCN grant DEC-2012/07/B/ST1/03417 
 (Poland).}\\
  \normalsize{Institute of Mathematics}\\
  \normalsize{University of Warsaw}\\
 \normalsize{ul.~Banacha 2, 02-097 Warszawa}\\
 \normalsize{Poland}\\
  \\
  {\normalsize{
  e-mail: annatal@mimuw.edu.pl}}
}
\date{}

\begin{document}
 \maketitle
  \begin{abstract}
We consider  standard $\La$-coalescents (or coalescents with multiple collisions) with a non-trivial ``Kingman part''. Equivalently, the driving measure $\Lambda$ has an atom at $0$; $\Lambda(\{0\})=c>0$. It is known that 
all such coalescents come down from infinity. 
Moreover, the  number of blocks $N_t$  is asymptotic to $v(t) = 2/(ct)$ as $t\to 0$. 
In the present paper we investigate  the second-order asymptotics of $N_t$ in the functional sense at small times. 
This  complements our earlier results on the fluctuations of the number of blocks for a class of regular $\La$-coalescents without the Kingman part.
In the present setting it turns out that 
the Kingman part dominates, and the limit process is a Gaussian diffusion, as opposed to the stable limit in our previous work. 
\\
\vglue .3cm
\noindent

 \end{abstract}
 \newpage
 
 \section{Introduction and main results}
 \subsection{Background}
  \label{S:i}
The  Kingman coalescent, introduced in \cite{Kingman1,Kingman2},  is one of the pillar processes of mathematical population genetics. 
The research reported here is linked to some of the classical results on the Kingman coalescent.
In particular, Griffiths in \cite{Griffiths1984} derives the Gaussian behavior of  the number of blocks (one dimensional distributions only).
Similar limits are discussed by Aldous in \cite{aldous_survey} in the absence of mutations, with general acknowledgement (as folk theorem), but no specific reference provided. More precisely, let $K_t$ be the number of blocks in the standard Kingman coalescent at time $t$. Then 
\cite{aldous_survey}  outlines the argument for 
\[
\sqrt{\frac{3t}{2}}\left( K_t - \frac{2}{t} \right) \Rightarrow N(0, 1),\ \mbox{ as } t\to 0.
\]

 The $\La$-coalescents form the simplest class of processes with exchangeable dynamics that generalize the Kingman coalescent.
They were introduced and first studied independently by Pitman \cite{pit99} and Sagitov \cite{sag99}, and were also considered in a contemporaneous work of Donnelly and Kurtz \cite{dk99}.
For recent overviews of the literature we refer the reader to \cite{bertoin, ensaios}. 

Let $\La$ be an arbitrary finite measure on $[0,1]$.
We denote by $(\Pi_t, t \ge 0)$ the associated $\Lambda$-coalescent. This Markov jump process $(\Pi_t, t \ge 0)$ takes values in the set of partitions of $\{1, 2, \ldots\}$. Its law is specified by the requirement that, for any $n\in \N$, the restriction $\Pi^n$ of $\Pi$ to $\{1, \ldots, n\}$ is a continuous-time Markov chain with 
the following transitions: whenever $\Pi^n$ has $ b \in \{2,\ldots,  n\}$ blocks, any given $k$-tuple of blocks coalesces at rate $\lambda_{b,k} := \int_{[0,1]} r^{k-2}(1-r)^{b-k} \La(dr)$.
The case $\La(dx) = \delta_0(dx)$ corresponds to the classical Kingman coalescent, where each pairwise collision occurs at rate 1, and no multiple collision is possible. 
The total mass of $\La$ can be scaled to $1$. 
This is convenient for the analysis, and corresponds to a constant time rescaling of the process. Henceforth we assume that $\La$ is a probability measure.
One of our main current assumptions is that 
$\La(\{0\}) =c >0$.
We distinguish two cases: if $c=1$ we call the corresponding coalescent the {\em pure Kingman} coalescent, while if 
$c\in (0,1)$ we call it the {\em mixed (with) Kingman} coalescent.

The {\em standard} $\La$-coalescent starts from the trivial configuration $\{\{i\}:i\in \N\}$.
We shall denote by $N_t$ 
the number of blocks of $\Pi(t)$ at time $t$. 
Note that the law of $N$ 
depends on $\La$, but it will be clear from the context which $\La$ (and therefore which $N$) we currently consider.
If $\P(N_t <\infty, \forall t > 0)=1$ the coalescent is said to \emph{come down from infinity} (CDI). 
It is well known that the Kingman coalescent has this property. Necessary and sufficient conditions for CDI for general $\Lambda$-coalescents were derived in \cite{sch1} and \cite{blg3}. In \cite{bbs} and \cite{bbl1} the small time behavior of $\Lambda$ coalescents was studied. In \cite{bbl1}, 
for a general $\Lambda$-coalescent that comes down from infinity, the authors found a non-random function $t\mapsto v_t$, dependent on $\Lambda$, such that as $t\to 0$
\begin{equation}\frac{N_t}{v_t}\to 1\qquad  \textrm{a.s.\ and in}\  L^p,
\label{e:speed}
\end{equation}
for any $p\ge 1$.  Any function $v$ satisfying \eqref{e:speed} is referred to as the \textit{speed of coming down from infinity}. There are many functions with this property, but clearly they have the same asymptotic behavior near $0$.
In our previous work \cite{LT} we investigated the second order asymptotics near $0$ for the number of blocks in a $\Lambda$-coalescent that comes down from infinity, assuming that $\Lambda$ has no atom at $0$ (no Kingman part) and that $\Lambda(\{1\})=0$. We studied the asymptotic behavior in a functional sense.
More precisely, we were interested in the processes
\begin{equation}
 r(\varepsilon)\left(\frac{N_{\varepsilon t}}{v_{\varepsilon t}}-1\right),\quad t\ge 0
\label{e:problem}
 \end{equation}
where $\varepsilon>0$, and $r(\varepsilon)$ is an appropriate norming, such that these processes converge in law in the Skorokhod space $D([0,\infty))$ as $\varepsilon\to 0$. 
We have shown that if $\Lambda$ is sufficiently regular near $0$, 
that is, if $\Lambda$ has a density in a neighborhood of $0$ that behaves as $Ay^{-\beta}$ with $0<\beta<1$ for some $A\in (0,\infty)$, then (for an appropriate speed  $v$) the correct norming is $r(\varepsilon)=\varepsilon^{-1/(1+\beta)}$.  Furthermore,  the limit process is a $(1+\beta)$-stable process of the form 
\begin{equation}K\frac{1}t{\int_{[0,t]}}udL_u,\label{e:stable}
\end{equation}
 where $L$ is a $(1+\beta)$-stable L\'evy process,  totally skewed to the left (it has no positive jumps),  and $K$ is a positive constant.

The object of the present paper is to present a complementary result, concerning the second order asymptotics of the number of blocks at small times for $\Lambda$-coalescents that have non-zero Kingman part: $\Lambda(\{0\})=c>0$. The presence of an atom at $0$ introduces some essential differences and the results of \cite{LT} cannot be applied to this case. However, as we will see, the main idea can be adapted to cover this case as well.

For other second-order fluctuation limits in the setting of exchangeable coalescents, we refer the reader to the works of Schweinsberg \cite{Schweinsberg2012},
Kersting \cite{Kersting},
Dahmer, Kersting and Wakolbinger \cite{DKW}, Kersting, Schweinsberg and Wakolbinger
\cite{KSW}.

\subsection{Main results}
Let $N=(N_t,\,t\ge 0)$ be the block counting process in a $\Lambda$-coalescent.
In the mixed with Kingman case, where $\La(\{0\}) =c >0$, it is easy to see  (by comparing with the Kingman coalescent slowed down by a factor $c$, which corresponds to the driving measure $\La(dx)=c\delta_0(dx)$) that  $\P(N_t<\infty) =1$ for all $t>0$, without any additional assumptions on $\Lambda$.
Furthermore, from the results of Berestycki et al.~\cite{bbl1} 
it follows that in this case \eqref{e:speed} is satisfied with the function $t\mapsto\frac{2}{ct}$, 
which therefore is 
 a {speed of CDI} for the corresponding $\La$-coalescent.
Note that this expression for the  speed depends only on the atom at $0$. In particular, the pure Kingman coalescent slowed down by a factor $c$ will have exactly the same speed of CDI.

The object of interest is the process \eqref{e:problem} as $\varepsilon\to 0$. We now set $r(\varepsilon)=\varepsilon^{-\frac 12}$.

We first consider the  pure Kingman case.
As already mentioned, a similar study had already been  undertaken in \cite{Griffiths1984, aldous_survey} in the setting of pure Kingman coalescent,  but only for the marginal distributions  (i.e. for fixed $t=1$, an analogue of the classical CLT).
In the present paper we study these fluctuations in a functional sense.
We then proceed to the general result for the mixed Kingman coalescent, which is novel even for the one dimensional distributions.

Let 
$D([0,\infty))$ denote the Skorokhod space of c\`adl\`ag functions equipped with the usual $J_1$ topology.

\medskip

Our main result in the pure Kingman setting is as follows:
\begin{thm}\label{thm:pureKingman}
Let $(N_t)_{t\ge 0}$ be the block counting process  in a standard Kingman coalescent. Then the process $X_{\varepsilon}$ defined by 
\begin{equation}
 X_{\varepsilon}(t)=\varepsilon^{-\frac 12}\left(\frac{\varepsilon t}{2}N_{\varepsilon t}-1\right), \quad t\ge 0, \qquad X_{\varepsilon}(0)=0
 \label{e:Xe}
\end{equation}
converges in law in $D([0,\infty))$ as $\varepsilon\to 0$ to a Gaussian process
\begin{equation}
 Z_t= \frac 1{\sqrt{2}\; t}\int_0^t udW_u,\quad t>0, \qquad Z_0=0, \label{e:limit_proc}
\end{equation}
where $W$ is a standard Brownian motion.
\end{thm}
\begin{rem}
 (a) The limit process $Z$ has the same form, as the one in \cite{LT}  (in the case where $c=0$, and  where $\Lambda$ has a density near zero, which behaves as $Cy^{-\beta}$; cf.~\eqref{e:stable}),
 if one formally sets
 $\beta=1$.
\\
(b) It is easy to see that the process $Z$ satisfies the equation
\begin{equation}
  Z_t=-\int_0^tZ_s\frac{1}{s} ds+\frac 1{\sqrt 2}W_t, \qquad Z_0=0.\label{e:eqZ}
\end{equation}
(c) It is worth  pointing out that the limit process \eqref{e:limit_proc}  also appeared in  the context of scaling limits related to hierarchical random walks (see \cite{BGT_number_variance}, Proposition 2.11).
\end{rem}

\medskip

To state our result in the general setting of 
 mixed Kingman coalescents  we first need to recall the speeds of CDI used in \cite{bbl1} and \cite{LT}, while introducing some additional notation.
 
Assume that $\Lambda(\{0\})=c>0$ (clearly $c\le 1$), hence $\Lambda$ has the form $\Lambda=c\delta_0+(1-c)\Lambda_1$, where for $0<c<1$, $\Lambda_1$ is the uniquely determined probability measure on $[0,1]$ such that $\Lambda_1(\{0\})=0$. For $c=1$ set $\Lambda_1\equiv 0$. 
Denote
\begin{equation}
\Psi(q)=\int_{[0,1]}(qy-1+(1-y)^q)\frac{\Lambda(dy)}{y^2}, \qquad q\geq 1,
\label{e:Psi}
\end{equation}
where  the function  $y\mapsto\frac {qy-1+(1-y)^q}{y^2}$ is  continuously extended on $[0,1]$, so that its value at  $y=0$ is  $\frac{q(q-1)}{2}$. 
Let  $\Psi_1$  denote the function given by \eqref{e:Psi}, with $\Lambda$ replaced by $\Lambda_1$.
In particular we have
\begin{equation}
 \Psi(q)=c\frac{q(q-1)}{2}+(1-c)\Psi_1(q), \qquad \forall q\geq 1.
 \label{e:Psi_Psi1}
\end{equation}

Similarly let 
\begin{equation}
\Psi^*(q)=\int_{[0,1]}(qy-1+e^{-qy})\frac{\Lambda(dy)}{y^2},\qquad q\geq 0.
\label{e:Psistar}
\end{equation}
These functions have already appeared in earlier papers (see e.g.~ \cite{bbl1} and \cite{LT} for some of their properties and a discussion on relation between $\Psi$ and $\Psi^*$). In particular, we know that the functions $\Psi$ and $q\mapsto \frac {\Psi(q)}{q}$ are increasing, and the same holds for  $\Psi^*$ in place of $\Psi$. 

From the assumption $c>0$ it follows that, for any $a>1$, the integral $\int_a^\infty\frac 1{\Psi(q)}dq$ is finite and the same is true for $\Psi^*$ (which is also a condition for CDI, see \cite{blg3}). As in \cite{LT}, we define the function $v: (0,\infty)\mapsto (0,\infty)$ by 
\begin{equation}
 \label{e:vt}
 t=\int_{v_t}^\infty\frac 1{\Psi(q)}dq.
\end{equation}
By \eqref{e:Psi_Psi1} and Lemma 2.1 in \cite{LT} it follows that $\int_1^\infty\frac 1{\Psi(q)}dq=\infty$, hence 
 $v_t> 1$ for all $t>0$. (Note that the assumption $\Lambda(\{1\})=0$ in the formulation of Lemma 2.1 in \cite{LT} was not used in the proof.)

Analogously to $v$, one can define $v^*$, by substituting  $\Psi^*$ for $\Psi$ in \eqref{e:vt}. 
From the results of \cite{bbl1} it follows that   $N_t/v_t^*\to 1$ as $t\to 0$, almost surely and in $L^p$ for any $p\ge 1$, hence $v_t^*$ is a speed of CDI. The same is true for $v$ in place of $v^*$. Moreover, both $v$ and $v^*$   are asymptotic to  $w$ near zero,  where 
\begin{equation*}
 w_t=\frac 2{ct},\qquad t\geq 0.
\end{equation*}

Let us denote
\begin{equation}
 X^v_{\varepsilon}(t)=\varepsilon^{-\frac 12}\left(\frac{N_{\varepsilon t}}{v_{\varepsilon t}}-1\right), \quad X^v_\varepsilon(0)=0.
\label{e:Xv}
 \end{equation}
 Similarly, let $X_\varepsilon^{v^*}$ and $X_\varepsilon^w$ be the processes defined as in \eqref{e:Xv} with $v$ replaced by $v^*$ and $w$, respectively.

The convergence result for the number of blocks of the mixed Kingman coalescent, when normalized by the speed $v$ or $v^*$, is analogous to Theorem \ref{thm:pureKingman}. The  only assumption made  on the measure $\Lambda$ is that $\Lambda(\{0\})>0$.  
However, if one wishes to replace the speed $v$ by the simpler function $w$, given above, then additional assumptions on the measure $\Lambda$ are necessary:

\begin{thm} \label{thm:mixedKingman} Assume that $\Lambda(\{0\})=c>0$ and write $\Lambda=c\delta_0+(1-c)\Lambda_1$, as above, $0<c\le 1$, $\Lambda_1\equiv 0$ if $c=1$.\\
(i) The processes $X^v_\varepsilon$ and $X^{v^*}_\varepsilon$ converge in law in $D([0,\infty))$ as $\varepsilon\to 0$ to the process $\sqrt{c}Z$, where $Z$ is defined in \eqref{e:limit_proc}.

\medskip 
\noindent (ii) Suppose additionally that the function $\Psi_{1}$ defined by \eqref{e:Psi} with $\Lambda$ replaced by $\Lambda_1$ satisfies
\begin{equation}
\lim_{q\to \infty}\frac{\Psi_1(q)}{q^{3/2}} =0.
\label{e:assumptionPsi}
\end{equation}
Then $X^w_\varepsilon$ converges in law in $D([0,\infty))$ as $\varepsilon\to 0$ to the process $\sqrt c Z$.
\end{thm}

\begin{rem}
\label{rem:afterThmMixed}
(a) If $c=1$, then part \textit{(ii)} clearly restates  Theorem \ref{thm:pureKingman}.
\\
(b) 
In part \textit{(i)} the measure $\Lambda_1$ can be completely arbitrary, the limit only depends on $\Lambda(\{0\})$, which shows that the Kingman part dominates.
\\
(c) Here we see the  same phenomenon as in \cite{LT}, that the speed of CDI has to be carefully chosen, and that we cannot always replace $v$ by $w$. 
In fact, condition \eqref{e:assumptionPsi} is sharp. This is discussed in more detail in  Remark 
\ref{rem:Robust} at the end of Section \ref{ssec:3.4} after the proof of part (ii) of the theorem.  

It is known (and easy to see) that  the asymptotic behavior of $\Psi_1(q)$ as  $q\to \infty$ depends quite strongly on the behavior of $\Lambda_1$ near $0$
(see for example \cite{LT}, Lemma 2.5). However, to ensure \eqref{e:assumptionPsi} one does not need to assume much about the regularity of $\Lambda_1$  near $0$. 
 \end{rem}

A simple sufficient condition for is now given. 
\begin{prop}\label{prop:Lambda}
 If   $\int_{[0,1]}\frac{1}{\sqrt y}\Lambda_1(dy)<\infty$, then \eqref{e:assumptionPsi} holds.
\end{prop}

Due to Remark \ref{rem:Robust}, it is easy to give examples of $\La_1$ where the hypothesis of the proposition is not satisfied and neither \eqref{e:assumptionPsi} nor the conclusion of Theorem \ref{rem:afterThmMixed} (ii) hold (e.g.~$\La_1$ is  Beta$(2-\alpha,\alpha)$ distribution for any   $\alpha\geq 3/2$).
The intuition is that for such $\La_1$ the mass near $0$ is not sufficiently strong (when compared to the atom at $0$) to change the class of speeds, but it is sufficiently strong to perturb the second-order asymptotics of $N_t$.

The proofs of Theorems \ref{thm:pureKingman} and \ref{thm:mixedKingman} use some of the main ideas of our recent paper \cite{LT}, where we studied the case 
$\Lambda(\{0\})=0$, and where $\Lambda$ was sufficiently regular at $0$. These techniques relied heavily on a representation of  $\Lambda$-coalescents satisfying $\Lambda(\{0\})=0$ via Poisson random measures. This representation has been already observed by Pitman
(see \cite{pit99}), but we reformulated it in \cite{LT} to better suit our needs.

At first sight the case of the (mixed) Kingman coalescent is different, since the same representation cannot be used.
 However, it turns out that if, somewhat artificially, one writes the effect of the Kingman part with the help of a different Poisson random measure, then many of the arguments used in \cite{LT} may be adapted to this case as well. 
 
 In particular, we begin by explicitly writing out  an integral equation for the number of blocks $N_t$. 
This equation 
involves an integral with respect to a certain Poisson random measure. In the setting where $\La(\{0\})\in(0,1)$,
this measure consists of two essentially different pieces: the first  corresponding  to the Kingman part (the atom at $0$), 
and the  second to  multiple collisions (the measure $\Lambda_1$). 
The latter piece, which accounts for the individual block coloring, 
was introduced and thoroughly studied in \cite{LT}.  We shall rely on the results of that analysis.

However, some of the technical estimates need to be done differently.
In a certain sense, the case $\Lambda(\{0\})>0$ is simpler, since the Kingman part dominates, and the limits are Gaussian. As in \cite{LT}, one has to consider terms resulting from the non-Kingman part, but now one can use less precise estimates of these terms. 
We also make use of a standard result (found e.g.~in \cite{ethierkurtz}), a version of a martingale central limit theorem, which ensures convergence in law of martingales  whose jumps are well controlled and whose skew brackets converge to a deterministic function.

\medskip
The remainder of the paper is organized as follows. Section \ref{sec:pure_kingman} contains the proof of  Theorem \ref{thm:pureKingman}. Section \ref{sec:mixed_Kingman} is dedicated to the proof of Theorem  \ref{thm:mixedKingman}, and it also contains a proof of Proposition \ref{prop:Lambda} and Remark \ref{rem:Robust}.

Throughout the paper $C, C_1, C_2, \ldots $  denote  positive constants, which may differ from  line to line.\\
The symbol $\Rightarrow$ denotes convergence in law in the Skorokhod space $D([0,\infty))$ equipped with $J_1$ topology.

\section{The pure Kingman case}
\label{sec:pure_kingman}
In this section we prove Theorem \ref{thm:pureKingman}. 

\medskip
Let us recall first the following deterministic easy lemma, 
that will be used frequently in the proofs.
\begin{lem}(\cite{bbl1}, Lemma 10)
\label{L:bbl1}
Suppose $f,g:[a,b]\mapsto{\mathbb R}$ are   c\`adl\`ag functions
such that
$\sup_{x\in[a,b]} \left|f(x) + \int_a^x
g(u)\,du\right| \leq K$,
for some $K<\infty$.
 If in addition
$f(x)g(x) > 0$, $x\in [a,b]$ whenever $f(x) \neq 0$,
 then
\[
\sup_{x\in[a,b]} \left|\int_a^x g(u)\,du\right|\leq K \ \mbox{
and } \sup_{x\in[a,b]} |f(x)| \leq 2K.
\]
\end{lem}

 \medskip
 Recall that now $\Lambda=\delta_0$, so that $N_t$ denotes  the number of blocks of the Kingman coalescent at time $t$.

The process $(N_t)_{t\ge 0}$ is a pure death continuous time Markov chain, and as such has a simple description.  
If it is at state $n$, $n\ge 2$, then it jumps to $n-1$ with intensity $\binom{n}{2}$. 
 For our purpose, it will be convenient to express this process with the help of a Poisson random measure. This  will facilitate the study of fine asymptotic behavior of $N$ near zero and it will allow us to  use some of the  standard techniques of the theory of 
 integration with respect to Poisson random measures. We refer to Chapter 8 of \cite{PZ} for a summary of the main properties of such integrals. 
We always take c\`adl\`ag versions of  martingales expressed as integrals with respect to a compensated Poisson random measure.

Denote $\Delta=\{(i,j)\in \Z_+^2: 1\le i<j\}$. We will often denote a typical element of $\Delta$ by $\kk$.
Let $\pi^0$ be a Poisson random measure on $\R_+\times \Delta$ with the intensity measure $\nu^0=\ell\otimes\sum_{(i,j)\in \Delta}\delta_{(i,j)}$, where $\ell$ is the Lebesgue measure on $\R_+$ and $\delta_{(i,j)}$ is the Dirac delta measure. 
In other words, 
$((\pi^0([0,t]\times \{\kk\}))_{t\ge 0})_{\kk\in \Delta}$ are independent Poisson processes with intensity $1$. 

The standard Kingman coalescent may be constructed from $\pi^0$ as follows:
Arrivals in the process indexed by $\kk=(i,j)$ correspond to potential times of coalescence of blocks currently labeled by $i$ and by $j$, but coalescence occurs only if there are at least $j$ blocks in the current configuration. 
More precisely, initially we have trivial configuration $\{\{1\},\{2\},...\}$ consisting of singleton blocks. After each coalescence event, the blocks are reordered according to  their smallest element. The $i$th and $j$th block in the current ordering coalesce into one block at the next arrival time of $\pi^0(\cdot\times \{(i,j)\})$. This construction is very much related to the Donnelly-Kurtz modified lookdown process, see \cite{dk99}.

 By $\hat \pi^0$ we denote the compensated Poisson random measure
\begin{equation*}
 \hat\pi^0 =\pi^0-\nu^0.
\end{equation*}
Let us also denote  $\Delta_k=\{(i,j)\in \Delta: 1\le i<j\le k\}$
for $k\in \Z_+$.  Note that $\# \Delta_k=\binom{k}{2}$.

\medskip

The following lemma is important for our analysis.
\begin{lem}\label{lem:2.2}
 Under the assumptions of Theorem \ref{thm:pureKingman} we have
 \begin{equation}
\frac {t}2N_t=1-\int_0^t \left(\frac s2 N_s-1\right)\frac 1s ds -M_t +R_t, \qquad t\ge 0,
\label{e:N}
 \end{equation}
 where 
 \begin{equation}
  M_t= \frac 12 \int_{[0,t]}\int_{\Delta}s\ind_{\Delta_{N_{s-}}}(\kk)\hat \pi^0(ds d\kk ), \qquad t\ge 0.
 \label{e:M}
 \end{equation}
 and  $R$ is a continuous process such that for any $T>0$ there exists $C>0$ such that
 \begin{equation} 
 E \sup_{s\le t}\abs{R_s}\le Ct\qquad {t\le T}.
 \label{e:R}
 \end{equation}
\end{lem}

\proof
As already mentioned,
the Kingman coalescent comes down from infinity, hence for any $0<r\le t$ we have $N_t\le N_r<\infty$, almost surely. 
We may and will assume that the coalescent is constructed using the procedure described before Lemma \ref{lem:2.2}.

Due to this construction, we have
$N_t = N_r - \int_{(r,t] \times \Delta} \ind_{\Delta_{N_{s-}}}(\kk) \pi^0(ds,d\kk)$.
Therefore, introducing the compensated $\pi^0$,
\begin{align}
 N_t
 =&N_r-\int_{(r,t]}\int_{\Delta}\ind_{\Delta_{N_{s-}}}(\kk)\nu^0(dsd\kk )-\int_{(r,t]}\int_{\Delta}\ind_{\Delta_{N_{s-}}}(\kk)\hat\pi^0(dsd\kk )\notag \\
 =&N_r-\int_r^t \frac{N_s(N_s-1)}{2}ds-\int_{(r,t]}\int_{\Delta}\ind_{\Delta_{N_{s-}}}(\kk)\hat\pi^0(dsd\kk ).
 \label{e:2.4a}
\end{align}
This is permissible, since the jumps of $N$ on $[r,t]$ are discrete (isolated). 
Clearly we have 
\begin{equation*}
 t N_t=rN_r+\int_r^t N_s\,ds+\int_r^t s\,dN_s.
\end{equation*}
Hence, using \eqref{e:2.4a} we obtain,
\begin{multline}
 \frac {t}2N_t= \frac {r}2N_r+\int_r^t \left(\frac {N_s}2 -s\frac{N_s(N_s-1)}{4}\right) ds\,
 - \int_{(r,t]}\int_{\Delta}\frac s2\ind_{\Delta_{N_{s-}}}(\kk)\hat \pi^0(ds d\kk ),
 \label{e:eqN}
\end{multline}
for any $ t\ge r$.
 If one formally plugs in $r=0$ in this final expression, one readily sees that the final term equals $M_t$ from (\ref{e:M}), and that the drift term can be written as 
\begin{equation*}
- \int_0^t \frac{sN_s}{2} \frac{1}{s}\left(\frac{sN_s}{2}-1 \right) \, ds + 
 \int_0^t \frac{s N_s}{4} \, ds.
\end{equation*}
The point is that for $s \approx 0$ we have $\frac{sN_s}{2}\approx 1$, and this explains the form of the drift in \eqref{e:N}, provided we can argue that the errors are small. 
We will in fact show that both $M$ and the integral in \eqref{e:N} are well defined, and that for any fixed $t>0$, as $r\to 0$, the left and the right hand side of \eqref{e:eqN} converge in probability to the corresponding left and right hand side of \eqref{e:N}. Due to the c\`adl\`ag property of all the processes under consideration, \eqref{e:N}  holds for all $t\geq 0$ simultaneously.

First we show that $M$ is well defined. 
Due to \cite{bbl1}, Theorem 2 we have
\begin{equation}
\lim_{t\to 0} E \sup_{s\le t}
\left(\frac{s}{2}N_s-1\right)^2=0.
\label{e:e1}
\end{equation}
 Hence
$\frac{r}{2}N_r\to 1$  in $L^2$ (this convergence also holds a.s.~and in any $L^p$, $p\ge 1$). 
Now \eqref{e:e1} and the fact that $N$ is non-increasing  immediately imply  that for any $T>0$
 \begin{equation}
  E\sup_{s\le T}\left(\frac s2 N_s\right)^2 <\infty.
  \label{e:e1a}
 \end{equation}

Using the definition of $\nu^0$ and \eqref{e:e1a} we hence obtain
\begin{align}
 E\int_{[0,t]}\int_{\Delta}s^2\ind_{\Delta_{N_s-}}(\kk)\nu^0(ds d\kk )
 =E \int_0^t s^2 \frac{N_{s}(N_s-1)}{2}ds
 \le& 
 Ct \quad\  \text{if } \
 t\le T .
 \label{e:e4}
\end{align}
Due to Theorem 8.23 in \cite{PZ} and \eqref{e:e4} we obtain that $M$ given by \eqref{e:M} is  a well defined square integrable martingale. Moreover, it satisfies $E M_t^2\le Ct$ for $t\le T$ 
Hence, by Doob's $L^2$ maximal
inequality
\begin{equation}
 E  \sup_{s\le t}M_s^2\le 4Ct \qquad for \ t\le T.
\label{e:e4a}
 \end{equation}

We observe that the last term on the right hand side of \eqref{e:eqN}  is equal to $M_t-M_r$, and from  \eqref{e:e4a} it follows that   $M_r\to 0$ in $L^2$ as $r\to 0$.

Let us now examine the drift term in  \eqref{e:eqN}. It can be written as
 \begin{align}
 A_r(t):=&\frac 12 \int_r^t N_sds -\frac 12\int_r^t s\frac{N_s(N_s-1)}{2} ds\notag\\
  =&-\int_r^t \frac{N_s}{2}\left(\frac s2 N_s-1\right)ds +\int_r^t\frac{s}{4}N_s ds.
 \label{e:e2}
  \end{align}
This allows us to improve \eqref{e:e1} in a similar way as it was done in \cite{LT}, Lemma 3.7 for  $\Lambda$-coalescents
without the Kingman part. More precisely,
using \eqref{e:eqN}, \eqref{e:e2} and Lemma \ref{L:bbl1} (with $g(s)=\frac{N_s}{2}(\frac{sN_s}{2}-1)$), for $r\le t\le T$  we have
\begin{equation*}
 \sup_{r\le s\le t}\abs{\frac{s}{2}N_s-1}
\le 2 \left(\abs{\frac{r}{2}N_r-1}+\abs{M_r}+\sup_{r\le s\le t}\abs{M_s}+\int_r^t\frac{s}{4}N_s ds \right).
\end{equation*}
Squaring both sides of the last expression, applying  expectation and
using  \eqref{e:e1}, \eqref{e:e1a} and \eqref{e:e4a}, 
 we obtain that for any $T>0$ there exists $C>0$ such that
 \begin{equation}
  E \sup_{s\le t}
\left(\frac{s}{2}N_s-1\right)^2\le C t\qquad t\le T.
\label{e:e3}
 \end{equation}
Estimate \eqref{e:e3}, 
  together with Jensen's inequality readily implies that the integral with respect to $ds$  in \eqref{e:N} is well defined    
for all $t$ simultaneously, almost surely.

Moreover,
 we can express the drift term $A_r$ of \eqref{e:e2} as
 \begin{equation*}
  A_r(t)=-\int_r^t \left(\frac s2 N_s-1\right)^2\frac 1s ds -\int_r^t \left(\frac s2 N_s-1\right)\frac 1sds +\frac 12 \int_r^t\frac{s}{2}N_s ds.
 \end{equation*}
 By \eqref{e:e3} and \eqref{e:e1a}
\begin{equation*}
 E \abs{A_r(t)+\int_r^t \left(\frac s2 N_s-1\right)\frac 1sds}\le C_1 t\qquad \text{for all }\ t\le T,
\end{equation*}
where $C_1$ does not depend on $r$.
 
This shows that, as $r\to 0$,  $A_r(t)$ converges in $L^1$ to $-\int_0^t \left(\frac s2 N_s-1\right)\frac 1sds+R_t$, where 
\begin{equation*}
 R_t=-\int_0^t \left(\frac s2 N_s-1\right)^2\frac 1s ds+\frac 12 \int_0^t\frac{s}{2}N_s ds.
\end{equation*}
Again \eqref{e:e3} and \eqref{e:e1a} yield \eqref{e:R}.
\qed

\bigskip

Recall \eqref{e:Xe},  let $M$ be the martingale defined by \eqref{e:M}, and define
 \begin{equation}
  Y_t=-\frac{1}{t}\int_{[0,t]} udM_u,\quad t>0, \qquad Y_0=0.
 \label{e:Yt1}
 \end{equation}
and
\begin{equation}Y_\varepsilon(t)=\varepsilon^{-\frac 12}Y(\varepsilon t).
 \label{e:Ye}
\end{equation}

\begin{lem}\label{lem:2.3}
The  process $(Y_t)_{t\in \R_+}$ satisfies the equation
\begin{equation}
 Y_t=-\int_0^tY_s\frac{1}{s} ds-M_t.
\label{e:Yt}
\end{equation}
Moreover, for any $T>0$ there exists $C>0$ such that for all $t\le T$ 
\begin{equation}
 E\sup_{s\le t}Y_s^2\le Ct
\label{e:Yta}
 \end{equation}
 and
\begin{equation}
\lim_{\varepsilon\to 0} E \sup_{t\le T}\abs{X_\varepsilon(t)-Y_{\varepsilon}(t)}=0.
\label{e:e5}
\end{equation}
\end{lem}

\proof
It is clear that $Y$ is well defined. Moreover, if we denote $H_t=tY_t$, then by the definition of $Y$ and $M$
we have
\begin{equation*}
H_t=-\int_{[0,t]} udM_u= -\frac 12 \int_{[0,t]}\int_{\Delta}s^2\ind_{\Delta_{N_{s-}}}(\kk)\hat \pi^0(ds d\kk ), \qquad t\ge 0.
\end{equation*}
Thus $H$ is a martingale with  quadratic  variation 
\begin{equation*}
 \left[H\right]_t=\frac 14 \int_0^t\int_{\Delta}s^4\ind_{\Delta_{N_{s-}}}(\kk)\pi^0(ds d\kk ), \qquad t\ge 0,
\end{equation*}
(cf.~Theorem 8.23 in \cite{PZ}).
Consequently, using \eqref{e:e1a} we obtain
\begin{align*}
 E\left[H\right]_t=&\frac 14 E\int_0^t\int_{\Delta}s^4\ind_{\Delta_{N_{s-}}}(\kk)\nu^0(ds d\kk )\\
 =&\frac 14E \int_0^t s^4 \frac{N_s(N_s-1)}{2}ds\\
 \le& C \int_0^t s^2 ds
 = \frac{C}{3} t^3.
 \end{align*}
In particular, 
\begin{equation}
 E Y_t^2=\frac{1}{t^2}E H_t^2=\frac{1}{t^2}E [H]_t\le \frac C{3} t.
\label{e:EYt2}
 \end{equation}
The identity \eqref{e:Yt} follows by simple integration by parts (note that $t\mapsto \frac{1}{t}$ is continuous and of finite variation on any interval $[a,b]$, $0<a<b$). The only subtle point 
is the lack of regularity of 
 $t\mapsto \frac 1t$  at $0$. This difficulty is easily overcome, by
writing first the formula for $Y_t-Y_r$, for any $0<r<t$, 
\begin{equation}
 Y_t-Y_r=\int_r^t\frac 1{s^2}\int_0^s udM_u ds-\int_{(r,t]}dM_s=-\int_r^t\frac 1s Y_sds-M_t+M_r, 
\label{e:2.19c}
 \end{equation}
and then letting $r\to 0$. 
Here we use \eqref{e:EYt2} to bound 
$\int_r^t\frac 1s \abs{Y_s}\,ds$ uniformly in $r>0$, implying that 
$\int_0^t\frac 1s Y_s\,ds$ exists in the absolute sense, almost surely.

Estimate \eqref{e:Yta} follows from \eqref{e:2.19c}, Lemma \ref{L:bbl1} and \eqref{e:e4a}.
\medskip

To prove \eqref{e:e5}  we set $X_t=\frac t2 N_t-1$ and 
observe that by
\eqref{e:N} and \eqref{e:Yt} we have
\begin{equation*}
 X_t-Y_t=-\int_0^t(X_s-Y_s)\frac 1sds +R_t.
\end{equation*}
Another application of Lemma \ref{L:bbl1} yields 
\begin{equation*}
 \sup_{s\le t}\abs{X_s-Y_s}\le 2\sup_{s\le t}\abs {R_s},
\end{equation*}
so \eqref{e:R} implies 
$E \sup_{t\le T}\abs{X_\varepsilon(t)-Y_{\varepsilon}(t)} \leq 2 C \sqrt{\varepsilon} T$, and hence \eqref{e:e5}.
 \qed
  
  \medskip
  
 We are now ready to proceed to the proof of the second order asymptotics of the number of blocks of the pure Kingman coalescent. 
 
 \proofof{of Theorem \ref{thm:pureKingman}}
 Due to \eqref{e:e5} and the symmetry of the law of $W$ it suffices to show that the process $-Y_\varepsilon$, 
given by \eqref{e:Yt1}--\eqref{e:Ye},  converges in law in $D([0,\infty))$ to the process $Z$ given by \eqref{e:limit_proc}.
 
We start by showing  that for $H_\varepsilon$,   where $H_\varepsilon(t) :=-tY_\varepsilon(t)=\frac{1}{\sqrt \varepsilon\; \varepsilon}\int_0^{t\varepsilon}udM_u$, $t\geq 0$ we have
 \begin{equation}
  \label{e:conv_tY}
  \left(H_\varepsilon(t)\right)_{t\ge 0}\Rightarrow \left( \frac 1{\sqrt 2}\int_0^t udW_u\right)_{t\ge 0}.
 \end{equation}
 For this we  use Theorem 1.4 in Chapter 7.1 of \cite{ethierkurtz}, as noted in the introduction.  
Observe  that 
$H_\varepsilon(t) = - {\varepsilon^{-\frac32}}H(\varepsilon t)$, 
where $H$ is taken from the proof of Lemma \ref{lem:2.3}.
Therefore  $H_\varepsilon(t)$ is again an $L^2$-martingale and it has the form
\begin{equation}
 H_\varepsilon(t)=\varepsilon^{-\frac 32}\frac 12\int_0^{\varepsilon t}\int_{\Delta} s^2 \ind_{\Delta_{N_{s-}}}(\kk)\hat\pi^0(dsd\kk ).
 \label{e:He}
\end{equation}
By the properties of the compensated Poisson integral we have
\begin{align}
 \<H_\varepsilon\>(t)=&\frac{1}{4\varepsilon^3}\int_0^{\varepsilon t}\int_{\Delta}s^4\ind_{\Delta_{N_{s-}}}(\kk)\nu^0(dsd\kk )\label{e:He1}\\
 =&\frac{1}{4\varepsilon^3}\int_0^{\varepsilon t}s^4\frac{N_s(N_s-1)}{2}ds\notag\\
 =&\frac 12\int_0^t s^2(\varepsilon s)^2\frac{N_{\varepsilon s}(N_{\varepsilon s}-1)}{4}ds.
 \label{e:skew_bracket}
\end{align}
We next verify the assumptions (b) of \cite{ethierkurtz}, Theorem 1.4 in Chapter 7.1,  with $c(t)=c_{11}(t)=\frac 12\int_0^tu^2du$, 
and $\<H_\varepsilon\>$ corresponding to $A^{(n)}$. 
Since $\<H_\varepsilon\>$ is continuous, we only need to prove that  $\<H_\varepsilon\>(t) $ converges  to $\frac 12\int_0^t u^2du$ in probability,  for each fixed $t>0$,
and that for any $T>0$
\begin{equation}
\lim_{\varepsilon\to 0} E\sup_{t\le T}\abs{H_\varepsilon(t)-H_\varepsilon(t-)}^2=0.
\label{e:jump}
\end{equation}
The first claim follows readily from \eqref{e:skew_bracket} and \eqref{e:e3}. 
Equality \eqref{e:jump}  is true  due to  \eqref{e:He}, since from  this representation of $H_\varepsilon$ it follows that  the jumps of $H_\varepsilon$ on $[0,T]$ are uniformly bounded by
$
\frac 12 \varepsilon^{-\frac{3}{2}}(\varepsilon T)^2.
$


This finishes the proof of \eqref{e:conv_tY}. To see that   convergence of $H_\varepsilon$ implies the required convergence of $-Y_\varepsilon$, one can apply the  argument from
\cite{LT}:
use the continuity of $t\mapsto 1/t$  away from $0$, and near $0$  use the  estimate \eqref{e:Yta} together with 
an analogous bound
$E\sup_{s\le t}\abs{Z_s}^2\le Ct$ for $t\le T$,
 where $Z$ is the limit process . Due to Lemma \ref{L:bbl1}
the latter bound follows from \eqref{e:eqZ} in the same way that \eqref{e:Yta} followed from \eqref{e:e4a} .
 See Steps 2-4 in the proof of Lemma 4.8 in \cite{LT} for more details. Note that here Step 3 simplifies due to \eqref{e:Yta}.\qed

\section{The mixed with Kingman case}
\label{sec:mixed_Kingman}
In this section we prove Theorem \ref{thm:mixedKingman}. 
In Section \ref{ssec:Xe} we  present an outline of the proof of Theorem \ref{thm:mixedKingman} for $X_\varepsilon^v$,
 in Section \ref{ssec:details} we prove the key technical lemmas needed to fill in this outline, and in Section \ref{ssec:3.3} we discuss the convergence of the processes $X_\varepsilon^{v^*}$ and $X_\varepsilon^w$. The final subsection contains the proof of Proposition \ref{prop:Lambda}.

\subsection{Outline of the proof of Theorem \ref{thm:mixedKingman} for $X_\varepsilon^v$.}
\label{ssec:Xe}
The proof combines the ideas from \cite{LT} and the proof of Theorem \ref{thm:pureKingman}, therefore we only briefly sketch it, omitting the details and concentrating on the differences. 
Proofs of the technical lemmas that require some new calculations  
(Lemmas \ref{lem:3.2}, \ref{lem:3.4} and \ref{lem:3.5}) will be given in  Section \ref{ssec:details}. 

We will use the Poissonian construction of $\Lambda$-coalescents (cf.~\cite{pit99}). 
More precisely, we make use of the notation introduced in the pure Kingman case, as well as the ``enriched'' Poisson random measure taken directly from \cite{LT}.
Recall \eqref{e:Psi}--\eqref{e:Psi_Psi1} and \eqref{e:Xv}. 
Suppose that we are given $\La$ as in the introduction,  and recall how $c$ and $\Lambda_1$ were  defined there.
 Let $\pic$ be a  Poisson random measure on $\R_+\times \Delta$ with intensity measure $c\nu^0$.
 It is defined analogously  to  $\pi^0$, introduced before Lemma \ref{lem:2.2}.   
 Let  $\piEc$ be a Poisson random measure on $\R_+\times[0,1]\times [0,1]^{\N}$ with intensity measure $(1-c)ds\frac{\Lambda(dy)}{y^2}\mu(d\x)$,
where $\x=(x_1,x_2,\ldots)$  and
       where $\mu$ is  a countable product of Lebesgue measures on  $[0,1]$ (or equivalently,  the law of an i.i.d.~sequence  of random variables, distributed uniformly on $[0,1]$) . 
It corresponds  to the non-Kingman part of the coalescent, and it is directly related to $\pi^E$ from \cite{LT}.    
 Since $\mu$ is a product of Lebesgue measures, we will usually abbreviate $\mu(d\x)$ as $d\x$. 
See \cite{LT} for more details and an interpretation.  
 
It is important to assume that $\pic$ and $\piEc$ are  independent. 
Then one can construct a version of the $\La$-coalescent by the following procedure  (blocks are again ordered according to their smallest element): 
(a) upon arrival of an atom $(t,\kk)$ of $\pic$, perform the collapsing of blocks as described above Lemma \ref{lem:2.2};
(b) upon arrival of an atom $(t,y,\x)$ of $\piEc$, the $j$-th block present in the configuration at time $t-$ is colored if and only if $x_j\leq y$. 
Once the colors are assigned, in order to form the configuration at time $t$, merge all the colored blocks into a single block, and leave the other (uncolored) blocks intact.

Following \cite{LT}, 
we define a function 
\begin{equation}
 f(k,y,\x)=\left(\sum_{j=1}^k\ind_{\{x_j\le y\}}-1\right)^+=\sum_{j=1}^k\ind_{\{x_j\le y\}}-1 +\ind_{\bigcap_{j=1}^k\{x_j>y\}}, 
\label{e:f}
 \end{equation}
which quantifies the decrease in the number of blocks during one coalescent event induced by $\piEc$, given that  $k$  blocks are present 
just before this event.

Observe that 
\begin{align}
 \Psi(k)=&c\frac{k(k-1)}{2}+(1-c)\Psi_1(k)
 \label{e:3.1}\\
 =&c\nu^0(\Delta_k)+(1-c)\int_{[0,1]}\int_{[0,1]^\N}f(k,y,\x)d\x\frac{\Lambda_1(dy)}{y^2},\label{e:3.1c}
\end{align}
since $\int_{[0,1]^\N}f(k,y,\x)d\x=E(\xi-1)^+$, where $\xi$ has the Binomial$(k, y)$ distribution.

From Corollary 15 in \cite{bbl1} the following lemma can be derived:
\begin{lem} \label{lem:Psi_v}
 \begin{align}
  \lim_{q\to \infty}\frac{\Psi_1(q)}{q^2}=&0,\label{e:3.14}
  \\
 \lim_{q\to \infty}\frac{\Psi(q)}{q^2}=\frac{c}{2},
\label{e:3.16}\\
  \lim_{t\to 0+} \frac {ct}2 v_t =1.\label{e:3.15}
 \end{align}
\end{lem}
More precisely, in \cite{bbl1} this was formulated for $\Psi^*$ and $v^*$, but their behavior 
at $\infty$ and $0+$ is the same as that of $\Psi$ and $v$, respectively. For completeness we include a short argument in 
 Section \ref{ssec:details}.

\medskip
We keep the convention that $\hat \pi$ denotes the compensated Poisson random measure $\pi$.
Using the Poissonian construction of the $\Lambda$-coalescent described above, then compensating, and applying \eqref{e:3.1c} we have
 \begin{align*}
N_t=&N_r-\int_r^t \Psi(N_s)\,ds - c \int_{(r,t] \times \Delta} \ind_{\Delta_{N_{s-}}}(\kk) \hpic(ds,d\kk)\\
& - (1-c)\int_{(r,t] \times [0,1] \times [0,1]^\N} f(N_{s-}, y, \x)
 \hpiEc(ds, dy, d\x), \ t\ge r.
 \end{align*}
Next, realizing that  \eqref{e:vt} implies $v_t'=-\Psi(v_t)$, for all $t>0$, one can obtain the following lemma
 in the same way as
  \eqref{e:eqN} or  \cite{LT}, Lemma 3.3. 
\begin{lem}\label{lem:3.1}
 For any $r>0$ and all $t\geq r$ we have
 \begin{align*}
  \frac{N_t}{v_t}=&\frac{N_r}{v_r}-\int_r^t\frac{N_s}{v_s}\left(\frac{\Psi(N_s)}{N_s}-\frac{\Psi(v_s)}{v_s}\right)ds\\
-&\int_{(r,t]}\int_\Delta \frac{\ind_{\Delta_{N_{s-}}}(\kk)}{v_s}\hpic(dsd\kk)
-\int_{(r,t]}\int_{[0,1]}\int_{[0,1]^\N}\frac{f(N_{s-},y,\x)}{v_s} \hpiEc(dsdyd\x). 
\end{align*}
\end{lem}

\medskip

As in the pure Kingman case, we wish to write the above equation starting from $r=0$. 
In particular, we need to show that
\begin{align}
 M_t^0:=&\int_{(0,t]}\int_\Delta \frac{\ind_{\Delta_{N_{s-}}}(\kk)}{v_s}\hpic(dsd\kk)\label{e:3.2}, \\
 M_t^1:=&\int_{(0,t]}\int_{[0,1]}\int_{[0,1]^\N}\frac{f(N_{s-},y,\x)}{v_s} \hpiEc(dsdyd\x), \label{e:3.3}
\end{align}
and
\begin{equation}
 A(t):=\int_0^t\frac{N_s}{v_s}\left(\frac{\Psi(N_s)}{N_s}-\frac{\Psi(v_s)}{v_s}\right)ds
\label{e:3.3a}
 \end{equation}
are all well defined. The integrals in \eqref{e:3.2} and \eqref{e:3.3} are to be understood in the sense of the usual (compensated) Poisson integration.

As before, from Theorem 2 in  \cite{bbl1} it follows that
\begin{equation}
 \lim_{t\to 0}E \sup_{s\le t}\left(\frac{N_s}{v_s}-1\right)^2=0
\label{e:3.5}.
 \end{equation}
We have already observed that $v_t\ge 1$ and $N_t$ is non-increasing, hence
 similarly to   \eqref{e:e1a} 
we find that 
there exists $C>0$ such that 
 \begin{equation}
 E\sup_{s\ge 0}\left(\frac{N_s}{v_s}\right)^2\le C.
\label{e:3.6}
 \end{equation}

In Section \ref{ssec:details} we will prove the following lemma.
\begin{lem}\label{lem:3.2} The processes 
 $M^0$ and $M^1$ given by \eqref{e:3.2} and \eqref{e:3.3} are well defined square integrable martingales. Moreover,  for any $T>0$ there exists $C>0$ such that 
 \begin{equation}
  E\sup_{s\le t} (M_s^i)^2\le Ct, \qquad t\le T, \quad i=0,1.
  \label{e:3.4}
 \end{equation}
\end{lem}

The function $q\mapsto \frac{\Psi(q)}{q}$ is increasing (see e.g.~Lemma 2.1(iv) in \cite{LT}), therefore, Lemma \ref{lem:3.1}, \eqref{e:3.6}, \eqref{e:3.4} and Lemma \ref{L:bbl1} imply
\begin{lem}\label{lem:3.3a} For any $T>0$ there exists $C>0$ such that 
\begin{equation}
 E\sup_{s\le t}\left(\frac{N_s}{v_s}-1\right)^2\le Ct, \qquad t\in [0,T].
\label{e:3.6a}
 \end{equation}
 \end{lem} 
The proof is very similar to that of \eqref{e:e3}.
 
 \medskip
 
Regarding the drift, we will prove the following in Section \ref{ssec:details}:
\begin{lem}
 \label{lem:3.4}
 For each $t>0$, the integral in \eqref{e:3.3a} is a well defined Lebesgue integral, almost surely.
 
 Moreover, 
 \begin{equation}
  A(t)=\int_0^t \left(\frac {N_s}{v_s}-1\right)\frac{1}{s}ds+U_t, \qquad t\ge 0,
 \label{e:3.6b}
 \end{equation}
where the process $U$ satisfies
\begin{equation}
 \label{e:3.8}
 \lim_{\varepsilon\to 0}\frac 1{\sqrt \varepsilon} E \sup_{s\le t}\abs{U_{\varepsilon s}}=0.
\end{equation}
\end{lem}

Note that this is simpler than the corresponding Lemma 4.4  in \cite{LT}, where at this point in the analogue of \eqref{e:3.6b}
we had to use a more complicated function instead of $\frac 1s$.

\medskip

We will also show in Section \ref{ssec:details} that the effect of $M^1$ is negligeable in the limit:
 \begin{lem}
  \label{lem:3.5}
  For any $1\le r<2$ and $T>0$ we have
  \begin{equation*}
   \lim_{\varepsilon\to 0}E \sup_{s\le T}\abs{\frac{1}{\sqrt \varepsilon}M^1_{\varepsilon t}}^r=0.
  \end{equation*}
 \end{lem}
 
 Using the above lemmas and Lemma \ref{L:bbl1} again, it is easy to deduce the following analogue of Lemma \ref{lem:2.2}:
\begin{lem}
 \begin{equation*} 
  \frac{N_t}{v_t}-1 =-\int_0^t \left( \frac{N_s}{v_s}-1\right)\frac 1sds -M_t^0 +R_t,
 \end{equation*}
 where,  for any $T>0$, $R$ satisfies
 \begin{equation}\label{e:3.10}
   \lim_{\varepsilon\to 0}E \sup_{s\le T}\abs{\frac{1}{\sqrt \varepsilon}R_{\varepsilon t}}=0.
\end{equation}
\end{lem}

The rest of the proof is the same as in the pure Kingman case.
Setting $R=0$ we study the process $(Y_t)$ satisfying 
\begin{equation*}
 Y_t =-\int_0^t Y_s \frac{1}{s}ds -M_t^0, \qquad Y_0=0,
\end{equation*}
which can be written explicitly as 
\begin{equation*}
 Y_t=-\frac{1}{t}\int_0^t sdM_s^0ds, \qquad t>0, \qquad Y_0=0.
\end{equation*}
From \eqref{e:3.10} and Lemma \ref{L:bbl1} it follows that the convergence in law in $D([0,\infty))$ of $X_\varepsilon^v$ is equivalent to the same convergence of 
\begin{equation*}
 Y_{\varepsilon}(t)=\frac{1}{\sqrt \varepsilon}Y_{\varepsilon t}.
\end{equation*}
This convergence is shown in exactly the same way as in the proof of Theorem \ref{thm:pureKingman}, using \eqref{e:3.5} and the fact that $v_t\sim \frac{2}{ct}$ as $t\to 0$.
 The constant $c$ appearing in the limit comes from the intensity of the Poisson random measure in the definition of $M^0$, 
which is $c\nu^0$ in this case,  compared to $\nu^0$ in the definition of $M$ in Section \ref{sec:pure_kingman}. In particular, the effect of this change is visible in  \eqref{e:He1}.

\subsection{Proofs of the auxiliary lemmas stated in Section \ref{ssec:Xe}}
\label{ssec:details}
\proofof {of Lemma \ref{lem:Psi_v}}
For $q\ge 1$ we can rewrite $\Psi_1(q)$ as
\begin{equation}
 \Psi_1(q)=q(q-1)\int_0^1\int_0^1\int_0^u(1-ry)^{q-2}drdu\Lambda_1(dy).
 \label{e:3.16b}
\end{equation}
Hence \eqref{e:3.14} follows since $(1-ry)^{q-2}\to 0$ as $q\to \infty$, and it is bounded by $1$ if $q\ge 2$.
Convergence \eqref{e:3.16} is a direct consequence of \eqref{e:3.1} and 
\eqref{e:3.14}. 

From the definition of $v_t$ we have that $ t v_t= v_t \int_{v_t}^\infty \frac{1}{\Psi(q)}dq $ , so that the l'Hospital rule and \eqref{e:3.16} imply 
\begin{equation*}
 \lim_{t\to 0+} \frac {ct}2 v_t
 =\frac c2\lim_{y\to \infty}y\int_{y}^\infty \frac{1}{\Psi(q)}dq
 =\frac c2\lim_{y\to \infty}\frac{-\frac{1}{\Psi(y)}}{-\frac{1}{y^2}}=1.
\end{equation*}
\qed

\proofof{of Lemma \ref{lem:3.2}}
 The argument for $M^0$ is the same as in Lemma \ref{lem:2.2}.
 Using \eqref{e:3.6}  we obtain the following analogue of \eqref{e:e4}:
\begin{equation}
 E \int_{(0,t]}\int_\Delta \left(\frac{\ind_{\Delta_{N_{s-}}}(\kk)}{v_s}\right)^2 c\nu^0(d\kk )ds
 =c\int_0^t E\left(\frac{N_s(N_s-1)}{2v_s^2}\right) ds\le C_1 t.
\label{e:3.7}
 \end{equation}
From the standard theory of Poisson integration (see e.g.~\cite{PZ}) it follows that $M^0$ 
is well defined. Moreover, it is a square integrable martingale with the following 
skew bracket 
\begin{align*}
 \<M^0\>_t=&c \int_0^t\int_\Delta \left(\frac{\ind_{\Delta_{N_{s-}}}(\kk)}{v_s}\right)^2 \nu^0(d\kk )ds
= c\int_0^t \frac{N_s(N_s-1)}{2v_s^2}ds.
\end{align*}
The bound \eqref{e:3.4} for $M^0$ now follows from Doob's $L^2$ maximal inequality and \eqref{e:3.7}.

\medskip

The argument for $M^1$ is similar. We need to get a bound on
\begin{equation*}
I(t):= E \int_0^t\int_{[0,1]}\int_{[0,1]^\N}\frac{f^2(N_{s-}, y,\x)}{v_s^2}d\x \frac{\Lambda_1(dy)}{y^2}dyds.
\end{equation*}
As $E \int_{[0,1]^\N}{f^2(k, y,\x)}d\x=E\left((\xi-1)^+\right)^2$, where $\xi$ is a Binomial$(k,y)$ r.v., it is elementary to check 
(see e.g~ Lemma 17 in \cite{bbl1} or (3.11) in \cite{LT}) that 
\begin{equation}
 \int_{[0,1]^\N}f^2(k,y,\x)d\x=k(k-1)y^2-k(k-1)\int_0^y\int_0^r(1-u)^{k-2}dudr.
 \label{e:3.13}
\end{equation}
By  
\eqref{e:3.13} and  \eqref{e:3.6}, for each $t>0$ we have
\begin{equation*}
I(t)\le E \int_0^t\int_{[0,1]}\frac{N_{s-}(N_{s-}-1)}{v_s^2}\Lambda_1(dy)ds
\le C t .
\end{equation*}
This implies that $M^1$ is well defined. Moreover, it is a square integrable martingale satisfying  $E( M_t^1)^2\le Ct$. 
As before, an application of Doob's $L^2$ maximal inequality finishes the proof.\qed

\medskip

\proofof{of Lemma \ref{lem:3.4}}
Fix any $T>0$, and suppose that  $t\in [0,T]$.
Observe that from \eqref{e:3.6a} and Jensen's inequality for $s\le t$ we have 
\begin{equation}E\abs{\frac{N_s}{v_s}-1}\le C_1\sqrt{s}.
\label{e:3.18} 
\end{equation}
Hence
\begin{equation}
E \int_0^t\abs{\frac {N_s-v_s}{v_s}}\frac 1sds\le C\sqrt t <\infty,
\label{e:3.18a}
\end{equation}
so the integral in \eqref{e:3.6b} is well defined.

Recalling \eqref{e:3.1}, and formally rewriting $A(t)$ defined by \eqref{e:3.3a}, we have
\begin{align}
 A(t)=&c\int_0^t\frac{N_s}{v_s}\frac{N_s-v_s}{2}ds+(1-c)\int_0^t\frac{N_s}{v_s}\left(\frac{\Psi_1(N_s)}{N_s}-\frac{\Psi_1(v_s)}{v_s}\right)ds\notag\\
 =&\int_0^t\frac{N_s-v_s}{v_s}\frac 1s ds+
  \int_0^t\frac{N_s-v_s}{v_s}\left(\frac {cs}2 v_s-1\right)\frac 1s ds \notag \\
 &+\frac c2\int_0^t \left(\frac{N_s-v_s}{v_s}\right)^2 v_s ds
 +(1-c) \int_0^t\frac{N_s}{v_s}\left(\frac{\Psi_1(N_s)}{N_s}-\frac{\Psi_1(v_s)}{v_s}\right)ds.
 \label{e:3.17}
\end{align}
It suffices to show that each of the terms in \eqref{e:3.17} is a well defined Lebesgue integral for all $t$ simultaneously, almost surely.
For this it is enough to show finiteness of $E\int_0^t\abs{\cdots}ds$ in each of the cases.  
The first term has  already been  estimated.

Considering the remaining terms in \eqref{e:3.17}, we denote
\begin{align}
 I_1(t)=&\int_0^t\abs{\frac{N_s-v_s}{v_s}(\frac {cs}2 v_s-1)}\frac 1s ds\label{e:3.19}\\
 I_2(t)=&\int_0^t \left(\frac{N_s-v_s}{v_s}\right)^2 v_s ds\label{e:3.20}\\
 I_3(t)=&\int_0^t\abs{\frac{N_s}{v_s}\left(\frac{\Psi_1(N_s)}{N_s}-\frac{\Psi_1(v_s)}{v_s}\right)}ds
 \label{e:3.21}
\end{align}
By \eqref{e:3.17}  and \eqref{e:3.18a}, the proof of the lemma will be completed once we 
show that
\begin{equation}
 E I_i(t)<\infty, \qquad i=1,2,3
 \label{e:3.22}
\end{equation}
and
\begin{equation}
\label{e:3.23}
\lim_{\varepsilon\to 0}\varepsilon^{-\frac 12} EI_i(\varepsilon T)=0,  \qquad i=1,2,3. 
\end{equation}
As already observed, $v_t$ is decreasing and $v_t\ge 1$, so from \eqref{e:3.15} it follows that there exists $C>0$ such that
\begin{equation}
 \label{e:3.24}
 v_s\le C\left(\frac 1s \vee 1 \right).
\end{equation}
In particular, this implies that $sv_s$ is bounded on $[0,T]$.

Estimates \eqref{e:3.24}, \eqref{e:3.6a} and \eqref{e:3.18} easily imply \eqref{e:3.22} for $i=1, 2$ and \eqref{e:3.23} for $i=2$. 
To show \eqref{e:3.23} for $i=1$ we additionally make an appropriate substitution and use the  dominated convergence theorem.

It remains to consider $I_3$. Let us denote $h_1(q)=\frac{\Psi_1(q)}{q}$. We can rewrite $h_1$ as
\begin{equation*}
 h_1(q)=\int_{[0,1]}\int_0^y (1-(1-r)^{q-1})dr\frac{\Lambda_1(dy)}{y^2}.
\end{equation*}
Using this representation, it is easy to see (cf.~\cite{LT}) that
\begin{equation}
\label{e:hprime}
\sup_{q>1} h_1'(q) \equiv \sup_{q>1} \abs{h_1'(q)} <\infty \mbox{ and } \lim_{q\to \iy} (h_1)'(q) = 0.
\end{equation}
By \eqref{e:3.21} and
the mean value theorem we have
\begin{equation*}
 EI_3(t)\le E\int_0^t\frac{N_s}{v_s}\abs{N_s-v_s}\sup_{q\ge N_s\wedge v_s}h_1'(q)ds.
\end{equation*}
Writing $\frac{N_s}{v_s}=(\frac{N_s}{v_s}-1)+1$, and estimating further,  we arrive at
\begin{equation*}
  EI_3(t)\le C E\int_0^t\left(\frac{N_s}{v_s}-1\right)^2v_sds+
  \int_0^t v_sE \left(\frac{\abs{N_s-v_s}}{v_s}\sup_{q\ge N_s\wedge v_s}h_1'(q)\right)ds.
\end{equation*}
Note that the first term on the right hand side is just $CI_2(t)$, which has  already been estimated. 
The second term is finite by \eqref{e:3.18} and \eqref{e:hprime}. 
To obtain \eqref{e:3.23} for $i=3$,  we apply the Cauchy-Schwarz inequality for the expectation inside the second integral, 
and then the dominated convergence theorem. 
Here we use \eqref{e:3.6a}, and the fact that
$E(\sup_{q\ge N_{\varepsilon s}\wedge v_{\varepsilon s}}h_1'(q))^2$ is bounded and tends to $0$ as $\varepsilon\to 0$, due to \eqref{e:hprime}, together with   $N_{\varepsilon s}\wedge v_{\varepsilon s}\to_{\varepsilon \to 0} \infty$, a.s.
This finishes the proof of \eqref{e:3.20} and \eqref{e:3.21} for $i=3$. \qed

\medskip

\proofof{of Lemma \ref{lem:3.5}}
Fix $T>0$, and let $0\le t\le T$.
We write 
\begin{equation}
 M_t^1=L_t+U_t,
 \label{e:3.30c}
\end{equation}
where
\begin{align*}
 L_t=& \int_0^t\int_{[0,1]}\int_{[0,1]^\N}y\hpiEc(dsdyd\x)\\
 U_t=&\int_0^t\int_{[0,1]}\int_{[0,1]^\N}\left(\frac
 {f(N_{s-},y,\x)}{v_s}-y\right)\hpiEc(dsdyd\x).
\end{align*}
Note that the process $(L_t)_{t\ge 0}$ is simply a L\'evy process with L\'evy measure $(1-c)\frac{\Lambda_1(dy)}{y^2}$. Both $U$ and $L$ are square integrable martingales.

By the standard properties of Poisson integrals (see e.g.~\cite{PZ}, Theorem 8.23) we have
\begin{equation}
 EU_t^2=
 E \int_0^t\int_{[0,1]}\int_{[0,1]^\N}\left(\frac
 {f(N_{s-},y,\x)}{v_s}-y\right)^2 d\x\frac{\Lambda_1(dy)}{y^2}ds
 \le 2 \left (J_1(t)+J_2(t)\right),
 \label{e:3.26a}
\end{equation}
where (using $N_s =N_ {s-} $ for almost all $s$) 
\begin{align*}
 J_1(t)=&E\int_0^t\int_{[0,1]}\int_{[0,1]^\N}
 \left(\frac{N_s}{v_s}\right)^2\left(\frac{f(N_{s},y,\x)}{N_s}-y\right)^2d\x \frac{\Lambda_1(dy)}{y^2}ds,  \mbox{and } \\
 J_2(t)=& E\int_0^t\int_{[0,1]}\left(\frac{N_s}{v_s}-1\right)^2 y^2 \frac{\Lambda_1(dy)}{y^2}ds.
\end{align*}
By \eqref{e:3.6a} we have 
\begin{equation}
 J_2(t)\le C t^2, \qquad \textrm{ for} \ t\le T. 
\label{e:3.27}
 \end{equation}
To estimate $J_1(t) $,  we recall  that for $k\ge 1$, $k\in \N$ 
\begin{equation*}
 \int_{[0,1]^\N} f(k,y,\x)d\x= k\int_0^y(1- (1-r)^{k-1})dr.
\end{equation*}
This, together with \eqref{e:3.13} gives
\begin{equation*}
\int_{[0,1]^\N} (f(k,y,\x)-ky)^2 d\x\le 2k^2y\int_0^y(1-r)^{k-1}dr.
\end{equation*}
Hence
\begin{equation*}
 J_1(t)\le 2E\int_0^t \int_{[0,1]}\left(\frac{N_s}{v_s}\right)^2\int_0^y(1-r)^{N_s-1}dr \frac{\Lambda_1(dy)}{y}ds.
\end{equation*}
 Due to \eqref{e:3.6} and  $(1-r)^{N_s-1}\le 1$, we get  $J_1(t)\le Ct$.

Moreover,
\begin{equation}
\label{e:Jone} 
 \frac{1}{\varepsilon}J_1(\varepsilon T)\le 2  \int_0^T\int_{[0,1]}\int_0^1E\left(\frac{N_{\varepsilon s}}{v_{\varepsilon s}}\right)^2(1-ry)^{N_{\varepsilon s}-1}dr\Lambda_1(dy)ds.
\end{equation}
Estimating $\left(\frac{N_{\varepsilon s}}{v_{\varepsilon s}}\right)^2\le 2\left(\frac{N_{\varepsilon s}}{v_{\varepsilon s}}-1\right)^2+2$,
 using $N_{\varepsilon s}\to \infty$ a.s.~and $0\le (1-ry)< 1$ a.e.~on the domain of integration, we have
\begin{equation*}
 E\left(\frac{N_{\varepsilon s}}{v_{\varepsilon s}}\right)^2(1-ry)^{N_{\varepsilon s}-1}
 \le 2 E\left(\frac{N_{\varepsilon s}}{v_{\varepsilon s}}-1\right)^2 +2E(1-ry)^{N_{\varepsilon s}-1}
 \to 0 \quad \textrm{as}\ \varepsilon\to 0.
\end{equation*}
Moreover, by \eqref{e:3.6}, the integrand in \eqref{e:Jone} is bounded for $s\le T$, $\varepsilon\le 1$, hence the dominated convergence theorem applies. 
We obtain that 
\begin{equation}
\label{e:3.28}
\lim_{\varepsilon\to 0} \frac{1}{\varepsilon}J_1(\varepsilon T)=0.
\end{equation}
By \eqref{e:3.26a}, \eqref{e:3.27} and \eqref{e:3.28}  we have
\begin{equation}
 \lim_{\varepsilon\to 0}E\left(\frac 1{\sqrt \varepsilon}U_{\varepsilon T}\right)^2=0.
\label{e:3.30a}
 \end{equation}
Since $U$ is a c\`adl\`ag martingale,  \eqref{e:3.30a} and Doob's $L^2$ maximal inequality imply that also
\begin{equation}
 \lim_{\varepsilon\to 0}E\sup_{t\le T}\left( \frac 1{\sqrt \varepsilon}U_{\varepsilon t}\right)^2=0.
\label{e:3.30b}
 \end{equation}
Let us now consider the process $L$.
Define
\begin{equation*}
 L_{\varepsilon}(t)=\frac{1}{\sqrt \varepsilon}L(\varepsilon t).
\end{equation*}
We will show that for any $0< r<2$ we have 
\begin{equation}
 \label{e:3.30}
 \lim_{\varepsilon\to 0}E \sup_{t\le T}\abs{L_{\varepsilon}(t)}^r=0.
\end{equation}
By Jensen's inequality it is clearly enough to consider $1<r<2$. Moreover, if $r\in(1,2)$ then the
Doob $L^r$ maximal inequality applied to the c\`adl\`ag martingale
 $L_\varepsilon$ implies that to obtain \eqref{e:3.30} it suffices to show that
\begin{equation}
 \label{e:3.31}
 \lim_{\varepsilon\to 0}E\abs{L_{\varepsilon}(T)}^r=0. 
\end{equation}

Fix any $r \in (1,2)$. In order to prove \eqref{e:3.31}, we will show that  $(\abs{L_\varepsilon(T)}^r)_{\varepsilon>0}$ is uniformly integrable, and that $L_\varepsilon(T)$ converges in law to $0$, and hence in probability.

The first assertion above follows from the simple estimate
\begin{equation*}
\sup_{\varepsilon> 0} E(L_{\varepsilon}(T))^2\le\sup_{\varepsilon>0} \frac{1}{\varepsilon}(1-c)\int_0^{\varepsilon T}\int_{[0,1]}y^2 \frac{\Lambda_1(dy)}{y^2}ds=(1-c)T.
\end{equation*}    

To show convergence in law of $L_\varepsilon(T)$, we write out its characteristic function:
\begin{align*}
 Ee^{iuL_{\varepsilon} (T)}=&\exp\left\{\varepsilon T\int_{[0,1]}\left( e^{iu\frac{1}{\sqrt \varepsilon}y}-1-iu\frac 1{\sqrt \varepsilon}y\right) \frac{\Lambda_1(dy)}{y^2}\right\}\\
 =&\exp\left\{ T u^2 \int_{[0,1]}\frac{ e^{iu\frac{1}{\sqrt \varepsilon}y}-1-iu\frac 1{\sqrt \varepsilon}y}{(\frac {y}{\sqrt \varepsilon})^2 u^2} {\Lambda_1(dy)} \right\}.
\end{align*}
Since $\sup_{x\in \R, x\neq 0} |\frac{e^{ix}-1-ix}{x^2}|$ is finite, and since 
$\lim_{ \abs{x}\to \infty} \frac{e^{ix}-1-ix}{x^2}=0$, the Lebesgue dominated convergence theorem implies that the right hand side above converges to $1$, 
as $\varepsilon\to 0$,  yielding the needed claim. 

Due to the previous discussion, \eqref{e:3.31} and therefore \eqref{e:3.30} holds .

The assertion of the  Lemma now follows by \eqref{e:3.30c}, \eqref{e:3.30b}, \eqref{e:3.30} and Jensen's inequality.
\qed

\medskip
\subsection{Convergence of $X_\varepsilon^{v^*}$.}
\label{ssec:3.3}

The proof is almost exactly the same as in the case of $\La$-coalescents without an atom at $0$,  studied in \cite{LT} (see Theorem 1.4 and Lemma 2.2 therein). 
Recalling the definitions of $X^{v^*}_\varepsilon$ and $X^v_\varepsilon$ we can write 
\begin{equation}
 X^{v^*}_\varepsilon(t)
 =X_\varepsilon^v(t)
 +\frac 1{\sqrt{\varepsilon}}\frac{N_{\varepsilon t}}{v_{\varepsilon t}}\left(\frac{v_{\varepsilon t}}{v^*_{\varepsilon t}}-1\right).
 \label{e:4.1}
\end{equation}
Therefore, to prove the desired convergence of $X^{v^*}_\varepsilon$, it suffices to show that
for any $T>0$ we have
\begin{equation}
\lim_{\varepsilon\to 0} E \sup_{t\le T}\frac 1{\sqrt{\varepsilon}}\frac{N_{\varepsilon t}}{v_{\varepsilon t}}\abs{\frac{v_{\varepsilon t}}{v^*_{\varepsilon t}}-1}=0.
\label{e:3.40}
 \end{equation}
From \eqref{e:3.6} it follows that \eqref{e:3.40} will hold  provided that
\begin{equation}
 \lim_{\varepsilon\to 0}\frac{1}{\sqrt\varepsilon}\left(\frac{v_{\varepsilon}}{v_\varepsilon^*}-1\right)=0.
\label{e:3.35a}
 \end{equation}
As in the proof of Theorem 1.4 in \cite{LT} (see Section 5.1 therein), the proof is thus 
reduced to a purely deterministic problem of showing \eqref{e:3.35a}.
 
Observe that 
\begin{equation*}
0\le  \Psi^*(q)-\Psi(q)=\frac{cq}2+(1-c)(\Psi_1^*(q)-\Psi_1(q)).
\end{equation*}
Lemma 2.1 in \cite{LT} implies that
\begin{equation}
\label{e:3.35}
 \Psi^*(q)-\Psi(q)\le \frac{q}{2}.
\end{equation}
This allows one to repeat the proof of Lemma 2.2 in \cite{LT}. \eqref{e:3.35} is used to obtain an estimate for the right hand side of   (2.15) in \cite{LT}.
 
From the analogue of (2.16) in \cite{LT} we obtain that there exists $t_0>0$ and $C>0$ such that for all $\varepsilon \in[0,t_0]$
\begin{equation*}
\abs{\frac{v_\varepsilon}{v_\varepsilon^*}-1}\le C\varepsilon.
\end{equation*}
Hence \eqref{e:3.35a} follows.

\subsection{Convergence of $X^w_\varepsilon$}\label{ssec:3.4}
Similarly to
   \eqref{e:4.1}, 
  we have 
\begin{equation}
 X^{w}_\varepsilon(t)
 =X_\varepsilon^v(t)
 +\frac 1{\sqrt{\varepsilon}}\frac{N_{\varepsilon t}}{v_{\varepsilon t}}\left(\frac{v_{\varepsilon t}}{w_{\varepsilon t}}-1\right).
 \label{e:4.5}
\end{equation}
As before, the proof of convergence of $X^w_\varepsilon$ reduces to showing that 
\begin{equation}
 \lim_{\varepsilon\to 0}\frac{1}{\sqrt\varepsilon}\left(\frac{v_{\varepsilon}}{w_\varepsilon}-1\right)=0.
\label{e:3.41}
 \end{equation}
This can be solved
using techniques from \cite{LT}, which rely on writing an equation for $\log\frac {v_\varepsilon}{w_\varepsilon}$ and applying Lemma \ref{L:bbl1} (cf. (2.15)-(2.16) in \cite{LT}).

Below we will use a different approach, which is   more  direct and at the same time shows that condition \eqref{e:assumptionPsi} is sharp.

We have the following lemma:
\begin{lem}\label{lem:4.1}
 Assume that
 \begin{equation}
\lim_{q\to \infty}\frac{\Psi_1(q)}{q^{3/2}} =A,
\label{e:assumptionPsiA}
\end{equation}
where $A\in [0,\infty]$. Then
 \begin{equation}
 \lim_{t\to 0+}\frac{1}{\sqrt t}\left(\frac {ct}2 {v_{t}}-1\right)=-\frac {2\sqrt 2}{3\sqrt  c}(1-c) A.
\label{e:4.2}
 \end{equation}
\end{lem}
\proof 
Using \eqref{e:vt} and an elementary integration we can write
\begin{equation}
\frac 1{\sqrt t}\left( v_t\frac{ct}{2}-1\right)=\frac c{2\sqrt t}v_t\left(t-\frac 2{c v_t}\right)
=(\frac c2)^{\frac 32} \frac 1{\sqrt {\frac c2 t v_t}} v_t^{\frac 32}\int_{v_t}^\infty\left( \frac 1{\Psi(q)}-\frac{1}{\frac c2 q^2} \right)dq.
\label{e:4.3}
\end{equation}
Applying  l'Hospital's rule, recalling \eqref{e:Psi_Psi1} and then using \eqref{e:3.16} and  \eqref{e:assumptionPsiA}, we obtain 
\begin{equation}
\lim_{y\to \infty}\frac {\int_{y}^\infty\left( \frac 1{\Psi(q)}-\frac{1}{\frac c2 q^2} \right)dq}{y^{-3/2}}
=\lim_{y\to \infty}\frac{\ \ \frac{\frac c2 y-(1-c)\Psi_1(y)}{\frac c2 y^2\Psi(y)}\ \ }{\frac 32 y^{-5/2}}
=-\frac 8{3c^2}A(1-c).
\label{e:4.4}
\end{equation}
Combining \eqref{e:3.15}, \eqref{e:4.3}  
and  \eqref{e:4.4} gives \eqref{e:4.2}.\qed

\medskip
An application of Lemma \ref{lem:4.1} with $A=0$ finishes the proof of convergence of $X_\varepsilon^w$.

\begin{rem}
\label{rem:Robust}
Observe that if the limit in \eqref{e:assumptionPsiA} exists and $A$ is  finite, then 
from \eqref{e:4.5}, Theorem \ref{thm:mixedKingman} (i), \eqref{e:3.5} and Lemma \ref{lem:4.1} it follows that 
$X_\varepsilon^w$ converges in law in $D([0,\infty))$, as $\varepsilon\to 0$, to the process 
\begin{equation*}
\tilde Z_t=\sqrt cZ_t-\frac{2\sqrt 2}{3\sqrt c}(1-c)A \sqrt t, \qquad t \ge 0,
\end{equation*}
where $Z$ is given by \eqref{e:limit_proc}. If \eqref{e:assumptionPsiA} holds with $A=\infty$ then $X_\varepsilon^w$ does not converge as $\varepsilon\to 0$.

It follows that the assumption \eqref{e:assumptionPsi} of Theorem \ref{thm:mixedKingman} (ii) may not be relaxed. In particular, 
in the Beta$(2-\alpha,\alpha)$ coalescent world, the parameter $\alpha=3/2$
is critical for the convergence in Theorem \ref{thm:mixedKingman} (ii).
\end{rem}

\subsection{Proof of Proposition \ref{prop:Lambda}}

Recall that $\Psi_1$ can be expressed as  in \eqref{e:3.16b}. Using the estimate $1-x\le e^{-x}$ and replacing the integral $\int_0^u\ldots dr$ by $\int_0^1\ldots dr$ one obtains
\begin{equation*}
 \Psi_1(q)\le q(q-1)\int_0^1\int_0^1 e^{-(q-2)ry}dr \Lambda_1(dy).
\end{equation*}
Moreover, observe that for any $\delta>0$ we have
\begin{align*}\sqrt q \int_0^1 \int_0^1 e^{-qry}dr &\Lambda_1(dy)
 \le \int_{[0,\delta]}\frac{1}{\sqrt y}\int_0^1\frac 1{\sqrt r}( \sqrt{qry}e^{-qry})dr \Lambda_1(dy)\\
 &+\int_{(\delta,1]}\frac{1}{y\sqrt{q}}\Lambda_1(dy)\\
 \le&  C\int_{[0,\delta]}\frac{1}{\sqrt y}\Lambda_1(dy)+\frac{1}{\delta \sqrt q},
\end{align*}
hence  the proposition follows.

\bigskip\noindent
\textbf{Acknowledgment.} Part of this work was done while A.\ Talarczyk was visiting Department of Mathematics, Universit\'e Paris-Sud. She is grateful for the hospitality of this institution.
We thank the reviewer and the associate editor for careful reading of the paper and helpful comments.

\end{document}